\documentclass[11pt]{article}
\usepackage{amssymb}
\usepackage{amsmath}
\usepackage{mathrsfs}
\usepackage{graphics}
\usepackage{graphicx}
\usepackage{xcolor}
\usepackage{subfigure}
\usepackage[T1]{fontenc}
\usepackage{latexsym,amssymb,amsmath,amsfonts,amsthm}\usepackage{txfonts}
\newcommand{\be}{\begin{eqnarray}}
\newcommand{\ben}{\begin{eqnarray*}}
\newcommand{\en}{\end{eqnarray}}
\newcommand{\enn}{\end{eqnarray*}}

\newtheorem{theorem}{Theorem}[section]
\newtheorem{lemma}{Lemma}[section]
\newtheorem{prp}[theorem]{Proposition}
\newtheorem{thm}[theorem]{Theorem}
\newtheorem{cor}[theorem]{Corollary}
\newtheorem{dfn}{Definition}[section]

\newtheorem{remark}{Remark}

\definecolor{rr}{rgb}{0,0,0}
\begin{document}
\renewcommand{\theequation}{\arabic{section}.\arabic{equation}}
\begin{titlepage}
\title{\bf Large deviations for quasilinear parabolic stochastic partial differential equations
}
\author{Zhao Dong$^{\dag,\ddag}$\ \ Rangrang Zhang$^{\flat}$,\ \ Tusheng Zhang$^{\sharp}$\\
{\small $^\dag$
RCSDS, Academy of Mathematics and Systems Science, Chinese Academy of Sciences, Beijing 100190, China.} \\
{\small $^\ddag$ School of Mathematical Sciences, University of Chinese Academy of Sciences.}\\
{\small $^\flat$ School of Mathematics and Statistics,
Beijing Institute of Technology, Beijing, 100081, China.}\\
{\small $^\sharp$ School of Mathematics, University of Manchester, Oxford Road, Manchester M13 9PL, England, UK}\\
({\sf dzhao@amt.ac.cn}, {\sf rrzhang@amss.ac.cn},\ {\sf tusheng.zhang@manchester.ac.uk} )}
\date{}
\end{titlepage}
\maketitle

\noindent\textbf{Abstract}:
In this paper, we establish the Freidlin-Wentzell's large deviations for quasilinear parabolic stochastic partial differential equations with multiplicative noise, which are neither monotone nor locally monotone. The proof is based on the weak convergence approach.



\section{Introduction}
In this paper, we are concerned with  quasilinear parabolic partial differential equations, which describe the phenomenon of convection-diffusion of the ideal fluids and therefore arise in a wide variety of important applications, including for instance two or three phase flows in porus media or sedimentation-consolidation processes (see, e.g. \cite{GMT} and the references therein ). The addition of a stochastic noise to this physical model is fully natural as it represents external random perturbations or a lack of knowledge of certain physical parameters. The
 quasilinear parabolic stochastic partial differential equations can be written as
\begin{eqnarray}\label{equ-0}
\left\{
  \begin{array}{ll}
   du+div(B(u))dt=div(A(u)\nabla u)dt+\sigma(u)dW(t),
  & \ x\in \mathbb{T}^d,\  t\in [0,T] \\
   u(0)=u_0. &
  \end{array}
\right.
\end{eqnarray}
Where $W(t), t\geq 0$ is a cylindrical Brownian motion, $A(\cdot), B(\cdot)$ are appropriate coefficients specified later. {\color{rr}The solution of (\ref{equ-0}) is denoted by $u=u(x,t)$.}
There are several recent works about the existence and uniqueness of pathwise weak solution of the above equation, i.e. strong in the probabilistic sense and weak in the PDE sense. We mention
two of them which are relevant to our work. In \cite{DHV}, Debussche, Hofmanov\'{a} and Vovelle obtained the well-posedness of the Cauchy problem for quasilinear degenerate parabolic stochastic partial differential equations.
    For the nondegenerate case, Hofmanov\'{a} and Zhang \cite{H-Z} provided a direct
    (and therefore much simpler) approach to the existence and uniqueness of the pathwise
     solution.

 The purpose of this paper is to establish the Freidlin-Wentzell's large deviation principle (LDP) for the quasilinear parabolic stochastic partial differential equations, which provides the exponential decay  of small probabilities associated  with
the corresponding  stochastic dynamical systems with small noise.

The proof of the large deviations will be based on the weak convergence approach introduced in Bou\'{e} and Dupuis \cite{MP}, Budhiraja and Dupuis \cite{BD}. As an important part of the proof, we need to obtain global well-posedness of the so called skeleton equation.
For the uniqueness, we adopt the method introduced in \cite{H-Z}. For the existence, we first use the heat kernel operators to construct a sequence of  approximating equations as in \cite{H-Z}. We then show that the family of the solutions of the approximating equations is compact in an appropriate space and that any limit of the approximating solutions gives rise to a solution of the skeleton equation. In this way, we remove the extra conditions imposed on the diffusion coefficient $\sigma(\cdot)$ in \cite{H-Z}. To complete the proof of the large deviation principle, we also need to study the weak convergence of the  perturbations of the
system (\ref{equ-0}) in the random directions of the Cameron-Martin space of the driving Brownian motions.

  \

This paper is organized as follows. The mathematical formulation of quasilinear parabolic stochastic partial differential equations is in Section 2. In Section 3, we introduce the  weak convergence method and state the main result. Section 4 is devoted to the study of the skeleton equations. The Large deviation principle is proved in Section 5.

\section{Framework}

{\color{rr} We work on a finite-time
interval $[0,T], T>0$ and consider periodic boundary conditions, that is, $x\in \mathbb{T}^d$ where $\mathbb{T}^d=[0,1]^d$ denotes the $d-$dimensional torus.}
Let $\mathcal{L}(K_1,K_2)$ (resp. $\mathcal{L}_2(K_1,K_2)$) be the space of bounded (resp. Hilbert-Schmidt) linear operators from a Hilbert space $K_1$ to another Hilbert space $K_2$, whose norm is denoted by $\|\cdot\|_{\mathcal{L}(K_1, K_2)}$(resp. $\|\cdot\|_{\mathcal{L}_2(K_1, K_2)})$. We will follow closely the framework of \cite{H-Z}. {\color{rr} In the following, we emphasize that the domain of all functional spaces is $\mathbb{T}^d$, for simplicity, we omit it. $C^1$ stands for the space of continuously differentiable functions. 
 $C^1_{lip}$ be the Lipschitz functions in $C^1$.} For $r\in[1,\infty]$, $L^r$ are the usual Lebesgue spaces and $\|\cdot\|_{L^r}$ represents the corresponding norm. In particular, denote $H=L^2(\mathbb{T}^d)$. For all $a\geq0$, let
$H^a$ be the usual Sobolev space of order $a$ with the norm
\[
\|v\|^2_{H^a}=\sum_{0\leq|\alpha|\leq a}\int_{\mathbb{T}^d}|D^{\alpha}v|^2dx.
\]
{\color{rr} Here, $\alpha$ is a multi-index, that is, $\alpha=(\alpha_1, \cdot\cdot\cdot, \alpha_d)$ with non-negative integers $\alpha_i, i=1,\cdot\cdot\cdot, d$. $|\alpha|=\sum^d_{i=1}\alpha_i$.}  $H^{-a}$ will denote the topological dual of $H^a$.

Now, we introduce the following hypothesis.
\begin{description}
   \item[\textbf{Hypothesis H}] \quad The flux function $B$, the diffusion matrix $A$, and the noise in (\ref{equ-0}) satisfy:
\begin{description}
  \item[(i)]
$B=(B_1,\cdot\cdot\cdot,B_d): \mathbb{R}\rightarrow \mathbb{R}^d$ is of class $C^1_{lip}$.

  \item[(ii)] $A=(A_{ij})^{d}_{i,j=1}: \mathbb{R}\rightarrow \mathbb{R}^{d\times d}$ is of class $C^1_{lip}$, uniformly positive definite and bounded, i.e. $\delta I\leq A\leq CI$.
  \item[(iii)] For each $u\in H$, $\sigma(u): U\rightarrow H$ defined by $\sigma(u)\bar{e}_k=\sigma_k(u(\cdot))$, where $U$ is a separable Hilbert space (with inner product $\langle\cdot,\cdot\rangle_U$ and norm $|\cdot|_U$), $(\bar{e}_k)_{k\geq 1}$ is an orthonormal basis of $U$ and  $\sigma_k(\cdot): \mathbb{R}\rightarrow \mathbb{R}$ are real-valued functions. In particular, assume that $\sigma$ satisfies the usual Lipschitz condition
       \begin{eqnarray}\label{equa-3}
      \sum_{k\geq 1}|\sigma_k(y)|^2\leq C(1+|y|^2),\quad  \sum_{k\geq1}|\sigma_k(y_1)-\sigma_k(y_2)|^2\leq C|y_1-y_2|^2 \quad {\rm{for}}\ y, y_1, y_2\in \mathbb{R}.
       \end{eqnarray}
\end{description}
\end{description}
Let $(\Omega, \mathcal{F}, \mathcal{F}_t, \mathbb{P})$ be a stochastic basis with a complete, right-continuous filtration with expectation $\mathbb{E}$. The driving process $W(t)$ is a $U-$cylindrical Wiener process defined on this stochastic basis whose paths belong to  $C([0,T],Y)$, where $Y$ is another Hilbert space such that the embedding $U\subset Y$ is Hilbert-Schmidt. $W$ admits the following decomposition
$W(t)=\sum^{\infty}_{k=1}\beta_k(t)\bar{e}_k$, $(\beta_k)_{k\geq 1}$ is a sequence of independent real-valued Brownian motions.
\begin{remark}
The above \textbf{(iii)} implies that $\sigma$ maps $H$ into $\mathcal{L}_2(U,H)$.
Thus, for a given predictable process $u\in L^2(\Omega, L^2([0,T],H))$, the stochastic integral $t\rightarrow \int^t_0\sigma(u(s))dW(s)$ is a well-defined $H-$valued square integrable martingale.
\end{remark}
Now, we recall the definition of a solution to (\ref{equ-0}) in \cite{H-Z}.
\begin{dfn}\label{dfn-1}
An $(\mathcal{F}_t)-$adapted, $H-$valued continuous process $(u(t), t\geq 0)$ is said to be a solution to equation (\ref{equ-0}) if
\begin{description}
  \item[(i)] $u\in L^2(\Omega, C([0,T],H)) \cap L^2(\Omega, L^2( [0,T], H^1))$,\ for any $T>0$,
  \item[(ii)] for any $\phi\in C^{\infty}(\mathbb{T}^d)$, $t>0$, the following holds almost surely
  \begin{eqnarray}\notag
  \langle u(t), \phi\rangle&-&\langle u_0, \phi\rangle-\int^t_0\langle B(u(s)), \nabla \phi\rangle ds\\
\label{equ-2}
  &=&-\int^t_0\langle A(u(s))\nabla u(s), \nabla \phi\rangle ds +\int^t_0\langle \sigma(u(s))dW(s),  \phi\rangle
  \end{eqnarray}
\end{description}
\end{dfn}
With the help of the global {\color{rr}well-posedness} results in \cite{H-Z}  and  a suitable approximation of initial values by smooth functions in \cite{DHV}, we have
\begin{thm}\label{thm-1}
Let $u_0\in L^p(\Omega,\mathcal{F}_0, L^p(\mathbb{T}^d))$ for all $p\in[1,\infty)$. Under the Hypothesis H, there exists a unique solution to the quasilinear SPDE (\ref{equ-0}) that satisfies the following energy inequality
\begin{eqnarray}\label{equ-3}
\mathbb{E}\sup_{0\leq t\leq T}\|u(t)\|^2_H+\int^T_0 \mathbb{E}\|u(t)\|^2_{H^1}dt<\infty.
\end{eqnarray}
\end{thm}
\begin{remark}
The solution here is pathwise, the so-called strong solution in the probabilistic sense.
\end{remark}
\section{The weak convergence approach and the statement of the main result}
In this section, we will recall the weak convergence approach developed  by Budhiraja and Dupuis in \cite{BD}. Let us first recall some standard definitions and results from the large deviation theory (see \cite{DZ})

Let $\{X^\varepsilon\}$ be a family of random variables defined on a probability space $(\Omega, \mathcal{F}, \mathbb{P})$ taking values in some Polish space $\mathcal{E}$.

\begin{dfn}
(Rate Function) A function $I: \mathcal{E}\rightarrow [0,\infty]$ is called a rate function if $I$ is lower semicontinuous. A rate function $I$ is called a good rate function if the level set $\{x\in \mathcal{E}: I(x)\leq M\}$ is compact for each $M<\infty$.
\end{dfn}
\begin{dfn}(Large Deviation Principle) The sequence $\{X^\varepsilon\}$ is said to satisfy a large deviation principle with rate function $I$ if for each Borel subset $A$ of $\mathcal{E}$
      \[
      -\inf_{x\in A^o}I(x)\leq \lim \inf_{\varepsilon\rightarrow 0}\varepsilon \log \mathbb{P}(X^\varepsilon\in A)\leq \lim \sup_{\varepsilon\rightarrow 0}\varepsilon \log \mathbb{P}(X^\varepsilon\in A)\leq -\inf_{x\in \bar{A}}I(x),
      \]
      where $A^o$ and $\bar{A}$ denote the interior and closure of $A$ in $\mathcal{E}$, respectively.
\end{dfn}

Suppose $W(t)$ is a $U$-cylindrical Wiener process defined on a probability space $(\Omega, \mathcal{F},\{\mathcal{F}_t\}_{t\in [0,T]}, \mathbb{P} )$. The paths of $W$ take values in $C([0,T],Y)$, where $Y$ is another Hilbert space such that the embedding $U\subset Y$ is Hilbert-Schmidt.

To state the criterion obtained by Budhiraja and Dupuis in \cite{BD}, we introduce the following spaces. Set
\begin{eqnarray*}
&\mathcal{A}=\Big\{\phi: \phi\ {\rm{is \ a \ U\text{-}valued }} \ \{\mathcal{F}_t\}-{\rm{predictable\ process\ such\ that}}\ \int^T_0 |\phi(s)|^2_Uds<\infty\ \mathbb{P}\text{-}a.s.\Big\};\\
&S_N=\Big\{ h\in L^2([0,T],U): \int^T_0 |h(s)|^2_Uds\leq N\Big\};\\
&\mathcal{A}_N=\{\phi\in \mathcal{A}: \phi(\omega)\in S_N,\ \mathbb{P}\text{-}a.s.\}.
\end{eqnarray*}
{\color{rr}Referring to \cite{BDM}, the set $S_N$ endowed with the following weak topology is a Polish space (complete separable metric space):
\begin{eqnarray*}
d_1(h,k)=\sum^{\infty}_{i=1}\frac{1}{2^i}\Big|\int^T_0\langle h(s)-k(s),\tilde{e}_i(s)\rangle_Uds\Big|,
\end{eqnarray*}
where $\{\tilde{e}_i(s)\}^{\infty}_{i=1}$ is a complete orthonormal basis of $L^2([0,T],U)$.} Then,
$S_N$ equippied with the above weak topology is a compact subspace of $L^2([0,T],U)$.
For $\varepsilon\geq 0$, suppose $\mathcal{G}^{\varepsilon}: C([0,T], Y)\rightarrow \mathcal{E}$ is a measurable map.  Set $X^{\varepsilon}=\mathcal{G}^{\varepsilon}(W)$ for $\varepsilon>0$. Consider the following conditions
\begin{description}
  \item[(I)] For every $N>0$, let $\{h_{\varepsilon}: \varepsilon>0\}$ $\subset \mathcal{A}_N$. If $h_{\varepsilon}$ converges to $h$ in distribution as $S_N-$valued random elements, then $\mathcal{G}^{\varepsilon}(W(\cdot)+\frac{1}{\sqrt{\varepsilon}}\int^{\cdot}_{0}h_\varepsilon(s)ds)$ converges to $\mathcal{G}^0(\int^{\cdot}_{0}h(s)ds)$ in distribution.
  \item[(II)] For every $N<\infty$, $\{\mathcal{G}^0(\int^{\cdot}_{0}h(s)ds): h\in S_N\}$ is a compact subset of $\mathcal{E}$.
\end{description}
The following result is due to Budhiraja and Dupuis in \cite{BD}.
\begin{thm}\label{thm-6}\cite{BD}
Suppose the above conditions (I), (II) hold. Then $X^{\varepsilon}$ satisfies a large deviation principle on $\mathcal{E}$ with a good rate function $I$ given by
\begin{eqnarray}\label{eq-5}
I(f)&=&\inf_{\{h\in L^2([0,T], U): f=\mathcal{G}^0(\int^{\cdot}_0 h(s)ds)\}}   \Big\{\frac{1}{2}\int^T_0|h(s)|^2_{U}ds\Big\},\ \ \forall f\in\mathcal{E}.
\end{eqnarray}
By convention, $I(\emptyset)=\infty.$
\end{thm}
Consider the following quasilinear parabolic stochastic partial differential equations driven by small multiplicative noise:
\begin{eqnarray}\label{equ-4}
\left\{
  \begin{array}{ll}
   du^{\varepsilon}(t)+div(B(u^{\varepsilon}(t)))dt=div(A(u^{\varepsilon}(t))\nabla u^{\varepsilon}(t))dt+\sqrt{\varepsilon}\sigma(u^{\varepsilon}(t))dW(t) , &\ x\in \mathbb{T}^d,\ t\in [0,T] , \\
    u^{\varepsilon}(0)=u_0.&
  \end{array}
\right.
\end{eqnarray}
 According to Theorem \ref{thm-1}, under the  Hypothesis H, there exists a unique strong solution of (\ref{equ-4}) in $C([0,T],H)\cap L^2([0,T],H^1)$. Therefore, there exists a Borel-measurable mapping
\[
\mathcal{G}^{\varepsilon}: C([0,T],Y)\rightarrow C([0,T],H)\cap L^2([0,T],H^1)
\]
such that $u^{\varepsilon}(\cdot)=\mathcal{G}^{\varepsilon}(W)$.

For $h\in L^2([0,T],U)$, consider the following skeleton equation
\begin{eqnarray}\label{equ-5}
\left\{
  \begin{array}{ll}
    du_h+div(B(u_h))dt=div(A(u_h)\nabla u_h)dt+\sigma(u_h)h(t)dt, &  \\
   u_h(0)=u_0. &
  \end{array}
\right.
\end{eqnarray}
The solution $u_h$, whose existence will be proved in next section,  defines a measurable mapping $\mathcal{G}^0: C([0,T],Y)\rightarrow C([0,T],H)\cap L^2([0,T],H^1)$ so that  $\mathcal{G}^0(\int^{\cdot}_0 h(s)ds):=u_h$.

Our main result reads as
\begin{thm}\label{thm-7}
Suppose the Hypothesis H is in place and $u_0 \in L^p(\mathbb{T}^d)$ for all $p\in[1,\infty)$. Then, $u^{\varepsilon}$ satisfies a large deviation principle on $C([0,T],H)$ with the good rate function $I$ defined by (\ref{eq-5}).
\end{thm}
 The rest of the paper is devoted to the proof of the main result.

\section{The skeleton equations}\label{s-1}
In this section, we will show that the skeleton equation (\ref{equ-5}) admits a unique solution for every $h\in L^2([0,T],U)$. The following result gives the uniqueness.
\begin{thm}\label{thm-8}
Assume the Hypothesis H holds. Then for every $h\in L^2([0,T],U)$, there exists at most one solution to equation (\ref{equ-5}) in the space $C([0,T],H)\cap L^2([0,T],H^1)$.
\end{thm}
\begin{proof} The proof is based on a suitable approximation of $L^1$ norm  introduced in \cite{H-Z}.
{\color{rr}
Let $1> a_1>a_2>\cdot\cdot\cdot>a_n>\cdot\cdot\cdot>0$ be a fixed sequence of decreasing positive numbers such that
\[
\int^{1}_{a_1}\frac{1}{r}dr=1,\ \cdot\cdot\cdot , \int^{a_{n-1}}_{a_n}\frac{1}{r}dr=n, \cdot\cdot\cdot
\]
Let $\psi_n(r)$ be a continuous function such that $supp(\psi_n)\subset (a_n, a_{n-1})$ and
\[
0\leq \psi_n(r)\leq 2\frac{1}{n}\times \frac{1}{r},\quad  \int^{a_{n-1}}_{a_n}\psi_n(r)dr=1.
\]
Define
\[
\phi_n(x)=\int^{|x|}_0\int^y_0\psi_n(r)drdy, \quad {\rm{for}}\ x\in \mathbb{R}.
\]
We have
\begin{eqnarray}\label{equ-00}
\phi_n(x)\leq |x|,  \ |\phi'_n(x)|\leq 1,\ 0\leq\phi''_n(x)\leq 2\frac{1}{n}\times \frac{1}{|x|},
\end{eqnarray}
and
\begin{eqnarray}\label{equ-46}
\phi_n(x)\rightarrow |x|, \quad {\rm{as}}\ n\rightarrow\infty.
\end{eqnarray}}
Define a functional $\Phi_n:H\rightarrow \mathbb{R}$ by
{\color{rr}
\[
\Phi_n(\gamma)=\int_{\mathbb{T}^d}\phi_n(\gamma(z))dz,\quad \gamma\in H.
\]
Then, we have
\begin{eqnarray*}
\Phi'_n(\gamma)(h)&=&\int_{\mathbb{T}^d}\phi'_n(\gamma(z))h(z)dz.
\end{eqnarray*}
}
Suppose that $u_1$, $u_2$ are two solutions to (\ref{equ-5}) with the same initial data $u_0$.
Applying the chain rule, we obtain
\begin{eqnarray*}
\Phi_n(u_1(t)-u_2(t))
&=& \int^t_0\Phi'_n(u_1(s)-u_2(s))d(u_1(s)-u_2(s))\\ \notag
&=& \int^t_0\int_{\mathbb{T}^d}\phi'_n(u_1(s,z)-u_2(s,z))\Big(-div(B(u_1(s,z)))+div(B(u_2(s,z)))\Big)dzds\\ \notag
&&+\int^t_0\int_{\mathbb{T}^d}\phi'_n(u_1(s,z)-u_2(s,z))\Big(div(A(u_1(s,z))\nabla u_1(s,z))-div(A(u_2(s,z))\nabla u_2(s,z))\Big)dzds\\ \notag
&&+\int^t_0\int_{\mathbb{T}^d}\phi'_n(u_1(s,z)-u_2(s,z))\Big(\sigma(u_1(s,z))h(s)-\sigma(u_2(s,z))h(s)\Big)dzds\\ \notag
&:=&I^1_n(t)+I^2_n(t)+I^3_n(t).
\end{eqnarray*}
Exactly argued as {\color{rr}(3.8) and (3.9)} in \cite{H-Z}, we have
\begin{eqnarray}\notag
I^1_n(t)+I^2_n(t)&\leq& \frac{C}{n}\int^t_0\int_{\mathbb{T}^d}|\nabla u_1(s,z)|dzds+\frac{C}{n}\int^t_0\int_{\mathbb{T}^d}|\nabla u_2(s,z)|dzds\\ \notag
&&+\frac{C}{n}\int^t_0\int_{\mathbb{T}^d}|\nabla u_2(s,z)|^2dzds+\frac{C}{n}\int^t_0\int_{\mathbb{T}^d}|\nabla u_1(s,z)|^2dzds\\
\label{equ-8}
&:=& K_n(t).
\end{eqnarray}
By H\"{o}lder's inequality, (\ref{equa-3}) and (\ref{equ-00}), we have
\begin{eqnarray}\notag
I^3_n(t)&\leq& \int^t_0\int_{\mathbb{T}^d}|\phi'_n(u_1(s,z)-u_2(s,z))||(\sigma(u_1(s,z))-\sigma(u_2(s,z)))h(s)|dzds\\ \notag
&\leq& \int^t_0\int_{\mathbb{T}^d}|\sum_{k\geq 1}\Big(\sigma_k(u_1(s,z))-\sigma_k(u_2(s,z))\Big)h_k(s)|dzds\\ \notag
&\leq& \int^t_0\int_{\mathbb{T}^d}(\sum_{k\geq 1}|\sigma_k(u_1(s,z))-\sigma_k(u_2(s,z))|^2)^{\frac{1}{2}}\Big(\sum_{k\geq 1}|h_k(s)|^2\Big)^{\frac{1}{2}}dzds\\ \notag
&\leq& C\int^t_0\int_{\mathbb{T}^d}|u_1(s,z)-u_2(s,z)||h(s)|_Udzds\\ \notag
&\leq& C\int^t_0|h(s)|_U\int_{\mathbb{T}^d}|u_1(s,z)-u_2(s,z)|dzds\\
\label{equ-9}
&=& C\int^t_0|h(s)|_U\|u_1(s)-u_2(s)\|_{L^1}ds.
\end{eqnarray}
Combining (\ref{equ-8}) and (\ref{equ-9}), we obtain

\begin{eqnarray}\label{equ-01}
\Phi_n(u_1(t)-u_2(t))
\leq{\color{rr} K_n(t)}+C\int^t_0|h(s)|_U\|u_1(s)-u_2(s)\|_{L^1}ds.
\end{eqnarray}
Let $n\rightarrow \infty$, we obtain
\begin{eqnarray*}
\|u_1(t)-u_2(t)\|_{L^1}
\leq C\int^t_0|h(s)|_U\|u_1(s)-u_2(s)\|_{L^1}ds.
\end{eqnarray*}
By Gronwall's inequality,  we get
\[
\|u_1(t)-u_2(t)\|_{L^1}=0,
\]
which implies the uniqueness.
\end{proof}

Next, we establish the existence of the solution of the skeleton equation (\ref{equ-5}).
\begin{thm}\label{thm-2}
Let $u_0\in H$. Under the Hypothesis H, for every $h\in L^2([0,T],U)$, there exists a solution to (\ref{equ-5}) that satisfies the following energy inequality
\begin{eqnarray}\label{equ-6}
\sup_{0\leq t\leq T}\|u_h(t)\|^2_H+\int^T_0\|u_h(t)\|^2_{H^1}dt<\infty.
\end{eqnarray}
\end{thm}
We will introduce some suitable approximating equations and show that the corresponding solutions converge to a solution to the skeleton equation. Let $P_r, r>0$ denote the semigroup on $H$ generated by the Laplacian on $\mathbb{T}^d$. Recall that
\[
P_r f(x)=\int_{\mathbb{T}^d}P_r(x,z)f(z)dz,
\]
where $P_r(x,z)$ stands for the heat kernel, $x, z\in\mathbb{T}^d$.

For $r>0$, $u\in H$, set
\begin{eqnarray}\label{equ-06}
A_r(u)(x)=P_r(A(u))(x), \ x\in\mathbb{T}^d,
\end{eqnarray}
where
\[
P_r(A(u))(x)=(P_r(A_{ij}(u))(x))^d_{i,j=1}.
\]
Note that there exists a constant $C$, independent of $r$,  such that
\begin{eqnarray}\label{equ-07}
\delta |\xi|^2\leq A_r(u)(x)\xi\cdot\xi\leq C|\xi|^2\quad \forall r>0,\  u\in H,\ x\in\mathbb{T}^d, \ \xi\in\mathbb{R}^d.
\end{eqnarray}
See (4.8) in \cite{H-Z}. Consider the following partial differential equations:
{\color{rr}\begin{eqnarray}\label{equ-08}
\left\{
  \begin{array}{ll}
   du^r_h(t)+div(B(u^r_h))dt=div(A_r(u^r_h(t))\nabla u^r_h(t))dt+\sigma(u^r_h(t))h(t)dt,  \\
   u^r_h(0)=u_0.
  \end{array}
\right.
\end{eqnarray}
We have the following result.
\begin{thm}\label{thm-3}
Let $u_0\in H$. Under Hypothesis H, for every $r>0$ and $h\in L^2([0,T],U)$, there exists a unique solution $u^r_h$ to equation (\ref{equ-08}). Moreover, the following uniform energy inequality holds
\begin{eqnarray}\label{equ-09}
\sup_{r}\Big\{\sup_{0\leq t\leq T}\|u^r_h(t)\|^2_H+\int^T_0\|u^r_h(t)\|^2_{H^1}dt\Big\}<\infty.
\end{eqnarray}
\end{thm}}
\begin{proof}
Set
\[
F_r(u):=-div(B(u))+div(A_r(u)\nabla u)+\sigma(u)h,\ u\in H^1.
\]
For $u\in H^1$, by H\"{o}lder's inequality and (\ref{equa-3}), we have
\begin{eqnarray*}
\langle\sigma(u(s))h(s), u(s)\rangle&=&\int_{\mathbb{T}^d}\sum_{k\geq 1}\sigma_k(u(s,z))h_k u(s,z)dz\\
&\leq& \int_{\mathbb{T}^d}\left(\sum_{k\geq 1}|\sigma_k(u(s,z))|^2\right)^{\frac{1}{2}}\left(\sum_{k\geq 1}|h_k|^2\right)^{\frac{1}{2}} |u(s,z)|dz\\
&\leq& |h(s)|_U\int_{\mathbb{T}^d}\left(\sum_{k\geq 1}|\sigma_k(u(s,z))|^2\right)^{\frac{1}{2}} |u(s,z)|dz\\
&\leq& |h(s)|_U\left(\int_{\mathbb{T}^d}\sum_{k\geq 1}|\sigma_k(u(s,z))|^2dz\right)^{\frac{1}{2}} \left(\int_{\mathbb{T}^d}|u(s,z)|^2dz\right)^{\frac{1}{2}} \\
&\leq& C|h(s)|_U\|u(s)\|_H+C|h(s)|_U\|u(s)\|^2_H,
\end{eqnarray*}
{\color{rr}where $h(s)=\sum^{\infty}_{k=1}h_k\bar{e}_k $ with $h_k=\langle h(s),\bar{e}_k \rangle_U$.}
Thus, by (4.10) in \cite{H-Z} {\color{rr}and using $x\leq 1+x^2$,} we have
\begin{eqnarray}\notag
\langle F_r(u(s)),u(s)\rangle&=&\langle B(u(s)), \nabla u(s)\rangle-\langle A_r(u(s))\nabla u(s), \nabla u(s)\rangle+\langle\sigma(u(s))h(s), u(s)\rangle\\ \notag
&\leq& C+C\|u(s)\|^2_H-\delta_1\|u(s)\|^2_{H^1}+C|h(s)|_U\|u(s)\|_H+C|h(s)|_U\|u(s)\|^2_H\\
\label{equ-20}
&\leq&C+C|h(s)|_U+C(1+|h(s)|_U)\|u(s)\|^2_H-\delta_1\|u(s)\|^2_{H^1},
\end{eqnarray}
for some constant $\delta_1>0$.
Moreover, by H\"{o}lder's inequality and (\ref{equa-3}),
\begin{eqnarray*}
&&\langle(\sigma(u(s))-\sigma(v(s)))h(s), u(s)-v(s)\rangle\\
&\leq& \int_{\mathbb{T}^d}|\sum_{k\geq 1}(\sigma_k(u(s,z))-\sigma_k(v(s,z)))h_k||u(s,z)-v(s,z)|dz\\
&\leq& \int_{\mathbb{T}^d}\left(\sum_{k\geq 1}|\sigma_k(u(s,z))-\sigma_k(v(s,z))|^2\right)^{\frac{1}{2}}\left(\sum_{k\geq 1}|h_k|^2\right)^{\frac{1}{2}}|u(s,z)-v(s,z)|dz\\
&\leq& |h(s)|_U\left(\int_{\mathbb{T}^d}\sum_{k\geq 1}|\sigma_k(u(s,z))-\sigma_k(v(s,z))|^2dz\right)^{\frac{1}{2}}\left(\int_{\mathbb{T}^d}|u(s,z)-v(s,z)|^2dz\right)^{\frac{1}{2}}\\
&\leq& C|h(s)|_U\|u(s)-v(s)\|^2_H,
\end{eqnarray*}
Combining the above estimates and recalling  (4.13) in \cite{H-Z}, we obtain
\begin{eqnarray}\notag
&&\langle F_r(u(s))-F_r(v(s)),u(s)-v(s)\rangle\\ \notag
&=&\langle B(u(s))-B(v(s)), \nabla (u(s)-v(s))\rangle-\langle A_r(u(s))\nabla u(s)-A_r(v(s))\nabla v(s), \nabla (u(s)-v(s))\rangle\\ \notag
&&\ +\langle(\sigma(u(s))-\sigma(v(s)))h, u(s)-v(s)\rangle\\ \notag
&\leq& C\|u(s)-v(s)\|^2_H+{\color{rr}C}\|u(s)-v(s)\|^2_H\|v(s)\|^2_{H^1}-\delta_2\|u(s)-v(s)\|^2_{H^1}+C|h(s)|_U\|u(s)-v(s)\|^2_H\\
\label{equ-21}
&\leq& C(1+|h(s)|_U+\|v(s)\|^2_{H^1})\|u(s)-v(s)\|^2_H-\delta_2\|u(s)-v(s)\|^2_{H^1},
\end{eqnarray}
for some constant $\delta_2>0$.
(\ref{equ-20}) and (\ref{equ-21}) together mean that the coefficient $F_r$ is locally monotone. Applying Theorem 1.1 in \cite{LR}, we obtain a unique solution $u^r_h$ to (\ref{equ-08}). By the chain rule,
\begin{eqnarray*}
\|u^r_h(t)\|^2_H&=&\|u_0\|^2_H+2\int^t_0\langle F_r(u^r_h(s)),u^r_h(s) \rangle ds\\
&\leq& \|u_0\|^2_H+C\int^t_0\Big(1+|h(s)|_U+(1+|h(s)|_U)\|u^r_h(s)\|^2_H\Big)ds
-\delta_1\int^t_0\|u^r_h(s)\|^2_{H^1}ds.
\end{eqnarray*}
Hence,
\begin{eqnarray}\notag
&&\|u^r_h(t)\|^2_H+\delta_1\int^t_0\|u^r_h(s)\|^2_{H^1}ds\\
\label{equuu-1}
&\leq&\|u_0\|^2_H+Ct+C\int^t_0|h(s)|_Uds+C\int^t_0(1+|h(s)|_U)\|u^r_h(s)\|^2_Hds.
\end{eqnarray}
By Gronwall's inequality, we obtain
\begin{eqnarray}\notag
&&\sup_{t\in [0,T]}\|u^r_h(t)\|^2_H+\delta_1\int^T_0\|u^r_h(s)\|^2_{H^1}ds\\
\label{equ-22}
&\leq & \Big(\|u_0\|^2_H+CT+C\int^T_0|h(s)|_Uds\Big)\cdot \exp\Big\{C\int^T_0(1+|h(s)|_U)ds\Big\}.
\end{eqnarray}
Because the constants involved in the above equation are independent of $r$, (\ref{equ-09}) follows.
\end{proof}

Let $K$ be a Banach space  with norm $\|\cdot\|_K$.
Given $p>1, \alpha\in (0,1)$, as in \cite{FG95}, let $W^{\alpha,p}([0,T]; K)$ be the Sobolev space of all $u\in L^p([0,T],K)$ such that
\[
\int^T_0\int^T_0\frac{\|u(t)-u(s)\|_K^{ p}}{|t-s|^{1+\alpha p}}dtds< \infty,
\]
endowed with the norm
\[
\|u\|^p_{W^{\alpha,p}([0,T]; K)}=\int^T_0\|u(t)\|_K^pdt+\int^T_0\int^T_0\frac{\|u(t)-u(s)\|_K^{ p}}{|t-s|^{1+\alpha p}}dtds.
\]
\begin{lemma}\label{lem-1}\cite{FG95}
Let $B_0\subset B\subset B_1$ be Banach spaces, $B_0$ and $B_1$ being reflexive, with compact embedding $B_0\subset B$. Let $p\in (1, \infty)$ and $\alpha \in (0, 1)$ be given. Let $X$ be the space
\[
X= L^p([0, T], B_0)\cap W^{\alpha, p}([0,T]; B_1),
\]
endowed with the natural norm. Then the embedding of $X$ into $L^p([0,T],B)$ is compact.
\end{lemma}
\begin{lemma}\label{lem-3}\cite{Temam-1}
 Let $V$ and $H$ be two Hilbert spaces  such that the imbedding  $V\subset H$ is compact. Denote by $V'$ the dual space of $V$. If $u\in L^2([0,T],V)$, $\frac{du}{dt}\in L^2([0,T],V')$, then $u\in C([0,T],H)$.
\end{lemma}
Next we prove the following result.
\begin{prp}\label{prpp-1}
Let $u^r$ be the solution to equation (\ref{equ-08}), then the family $\{u^r\}_{r\geq 0}$ is compact in $L^2([0,T],H)$.
\end{prp}
\begin{proof}
Recall
\begin{eqnarray*}
 u^r(t)&=&u_0-\int^t_0div(B(u^r(s)))ds+\int^t_0div(A_r(u^r(s))\nabla u^r(s))ds+\int^t_0\sigma(u^r(s))h(s)ds\\
&:=& I_1+I_2(t)+I_3(t)+I_4(t).
\end{eqnarray*}
Clearly, $\|I_1\|_H\leq C_1$. By integration by parts and Hypothesis H, we have
\begin{eqnarray*}
\|div(B(u^r(s)))\|_{H^{-1}}&=& \sup_{\|v\|_{H^1}\leq 1}|\langle v,div(B(u^r(s)))\rangle|\\
&=& \sup_{\|v\|_{H^1}\leq 1}|\langle \nabla v, B(u^r(s))\rangle|\\
&\leq& C(1+\|u^r(s)\|_H).
\end{eqnarray*}
Then
\begin{eqnarray*}
\|I_2(t)-I_2(s)\|^2_{H^{-1}}
&=&\|\int^t_sdiv(B(u^r(l)))dl \|^2_{H^{-1}}\\
&\leq& C(t-s)\int^t_s\|div(B(u^r(l)))\|^2_{H^{-1}}dl\\
&\leq& C(t-s)^2(1+\sup_{t\in [0,T]}\|u^r(t)\|^2_H),
\end{eqnarray*}
Hence, by Theorem \ref{thm-3}, we have for $\alpha\in (0,\frac{1}{2})$,
\begin{eqnarray*}\notag
&&\|I_2\|^2_{W^{\alpha,2}([0,T];H^{-1})}\\ \notag
&\leq&\int^T_0 \|I_2(t)\|^2_{H^{-1}}dt+\int^T_0\int^T_0\frac{\|I_2(t)-I_2(s)\|^2_{H^{-1}}}{|t-s|^{1+2\alpha}}dsdt\\
\label{equu-2}
&\leq& C_2(\alpha).
\end{eqnarray*}
Moreover, by the boundedness of $A_r$,
\begin{eqnarray*}
\|div(A_r(u^r(s))\nabla u^r(s))\|_{H^{-1}}&=& \sup_{\|v\|_{H^1}\leq 1}|\langle v,div(A_r(u^r(s))\nabla u^r(s))\rangle|\\
&=& \sup_{\|v\|_{H^1}\leq 1}|\langle \nabla v, A_r(u^r(s))\nabla u^r(s)\rangle|\\
&\leq& C(1+\|u^r(s)\|_{H^1}),
\end{eqnarray*}
and hence,
\begin{eqnarray*}
\|I_3(t)-I_3(s)\|^2_{H^{-1}}
&=&\|\int^t_s div(A_r(u^r(l))\nabla u^r(l))dl \|^2_{H^{-1}}\\
&\leq & C(t-s)\int^t_s\| div(A_r(u^r(l))\nabla u^r(l))\|^2_{H^{-1}}dl\\
&\leq& C(t-s)(T+\int^T_0\|u^r(s)\|^2_{H^1}ds).
\end{eqnarray*}
this again implies
\begin{eqnarray*}\label{equu-1}
\|I_3\|^2_{W^{\alpha,2}([0,T];H^{-1})}\leq C_3(\alpha).
\end{eqnarray*}
By Hypothesis H, we have
\begin{eqnarray*}
\|\sigma(u^r(s))h(s)\|_{H}&\leq& C|h(s)|(1+\|u^r(s)\|_{H}).
\end{eqnarray*}
Hence,
\begin{eqnarray*}
\|I_4(t)-I_4(s)\|^2_{H}
&=&\|\int^t_s \sigma(u^r(l))h(l)dl \|^2_{H}\\
&\leq & C(t-s)\int^t_s\| \sigma(u^r(l))h(l)\|^2_{H}dl\\
&\leq& C(t-s)(1+\sup_{s\in[0,T]}\|u^r(s)\|^2_{H})\int^T_0|h(s)|^2ds,
\end{eqnarray*}
which yields that for $\alpha\in (0,\frac{1}{2})$,
\begin{eqnarray*}\label{equu-3}
\|I_4\|^2_{W^{\alpha,2}([0,T];H)}\leq C_4(\alpha).
\end{eqnarray*}
Combining the above estimates together, we have for $\alpha\in (0,\frac{1}{2})$,
\begin{eqnarray}\label{equu-4}
\sup_r\|u^r\|^2_{W^{\alpha,2}([0,T];H^{-1})}< +\infty.
\end{eqnarray}
 By Theorem \ref{thm-3} and (\ref{equu-4}), we see that  $u^r$, $r>0$  are bounded uniformly in the space
\[
L^2([0,T],H^1)\cap W^{\alpha,2}([0,T];H^{-1}).
\]
Apply Lemma \ref{lem-1} to conclude the proof.
\end{proof}

From Theorem \ref{thm-3} and the proof of Proposition \ref{prpp-1} we see that the family $\{u^{r}\}_{r>0}$ is bounded in the space $C([0,T],H)\cap W^{\alpha,2}([0,T];H^{-1})$. Therefore, the following result is valid.
\begin{cor}\label{corr-1}
The family $\{u^{r}\}_{r>0}$ is tight in $C([0,T],H^{-1})$.
\end{cor}

As a consequence of Theorem \ref{thm-3}, Proposition \ref{prpp-1} and Corollary \ref{corr-1}, we have
\begin{cor}\label{cor-4}
There exist a sequence $\{u^{r_n}\}_{n\geq 1}$ and an element $u\in L^{\infty}([0,T],H)\cap L^2([0,T],H)\cap L^2([0,T],H^1)\cap C([0,T],H^{-1})$ such that
\begin{eqnarray*}
&u^{r_n}&\rightarrow u \ in\ L^{\infty}([0,T],H) \ \mbox{ {\color{rr}in weak star topology}},\\
&u^{r_n}&\rightarrow u \ strongly \ in\ L^2([0,T],H),\\
&u^{r_n}&\rightarrow u \ strongly \ in\ C([0,T],H^{-1}),\\
&u^{r_n}&\rightharpoonup u \ weakly \ in\ L^2([0,T],H^1).
\end{eqnarray*}
\end{cor}
Now, we are ready to prove Theorem \ref{thm-2}.

\begin{flushleft}
\textbf{Proof of Theorem \ref{thm-2}}. \quad Let $u$ be defined as in Corollary \ref{cor-4}. We will show that $u$ is a solution to the skeleton equation (\ref{equ-5}) and $u\in C([0,T],H)$.
For  a test function $\phi\in C^{\infty}(\mathbb{T}^d)$, we have
\begin{eqnarray}\notag
&&\langle u^{r_n}(t), \phi \rangle-\langle u_0, \phi \rangle-\int^t_0\langle B(u^{r_n}(s)), \nabla \phi \rangle ds\\
\label{equua-1}
&=&-\int^t_0\langle A_{r_n}(u^{r_n}(s))\nabla u^{r_n}(s), \nabla \phi \rangle ds+\int^t_0 \langle\sigma(u^{r_n}(s))h(s), \phi\rangle ds.
\end{eqnarray}
{\color{rr}Taking Corollary \ref{cor-4} into account}, we deduce that
\[
|\langle u^{r_n}(t)-u(t), \phi \rangle|\rightarrow 0.
\]
Since $u^{r_n}\rightarrow u \ strongly \ in\ L^2([0,T],H)$, we have
\[
|\int^t_0\langle B(u^{r_n}(s))-B(u(s)), \nabla \phi \rangle ds|\leq \|\nabla \phi\|_{L^{\infty}}\int^t_0\|u^{r_n}(s)-u(s)\|_Hds\rightarrow 0,
\]
and
\[
|\int^t_0 \langle(\sigma(u^{r_n}(s))-\sigma(u(s)))h(s), \phi\rangle ds|
\leq \| \phi\|_H\left(\int^t_0|h(s)|^2_Uds\right)^{\frac{1}{2}}\left(\int^t_0\|u^{r_n}(s)-u(s)\|^2_Hds\right)^{\frac{1}{2}}\rightarrow 0.
\]
Note that
\begin{eqnarray*}\notag
&&\int^t_0\langle A_{r_n}(u^{r_n}(s))\nabla u^{r_n}(s)-A(u(s))\nabla u(s), \nabla \phi \rangle ds\\
&=& \int^t_0\langle (A_{r_n}(u^{r_n}(s))-A_{r_n}(u(s)))\nabla u^{r_n}(s), \nabla \phi \rangle ds
\\
&&\ +\int^t_0\langle (A_{r_n}(u(s))-A(u(s)))\nabla u^{r_n}(s), \nabla \phi \rangle ds\\
&&\ +\int^t_0\langle A(u(s))(\nabla u^{r_n}(s)-\nabla u(s)), \nabla \phi \rangle ds.
\end{eqnarray*}
By the contraction property of the semigroup of $P_{r_n}$, Lipschitz continuity of $A$, Corollary \ref{cor-4} and Theorem \ref{thm-3}, we have
\begin{eqnarray*}
&&|\int^t_0\langle (A_{r_n}(u^{r_n}(s))-A_{r_n}(u(s)))\nabla u^{r_n}(s), \nabla \phi \rangle ds|\\
&\leq& \|\nabla\phi\|_{L^{\infty}}\int^t_0\int_{\mathbb{T}^d}|P_{r_n}(A(u^{r_n}(s))-A(u(s)))(x)||\nabla u^{r_n}(s)(x) |dxds\\
&\leq&
C\left(\int^t_0\|A(u^{r_n}(s))-A(u(s))\|^2_Hds\right)^{\frac{1}{2}}\left(\int^t_0\|u^{r_n}(s)\|^2_{H^1}ds\right)^{\frac{1}{2}}\\
&\leq&
C\left(\int^t_0\|u^{r_n}(s)-u(s)\|^2_Hds\right)^{\frac{1}{2}}\left(\int^t_0\|u^{r_n}(s)\|^2_{H^1}ds\right)^{\frac{1}{2}}\rightarrow 0.
\end{eqnarray*}
By the strong continuity of the semigroup $P_{r_n}$ and the boundedness of $A(u)$, we have
\begin{eqnarray*}
&&|\int^t_0\langle (A_{r_n}(u(s))-A(u(s)))\nabla u^{r_n}(s), \nabla \phi \rangle ds|\\
&\leq& \|\nabla\phi\|_{L^{\infty}}\int^t_0\int_{\mathbb{T}^d}|P_{r_n}A(u(s))(x)-A(u(s))(x)||\nabla u^{r_n}(s)(x) |dxds\\
&\leq&
C\left(\int^t_0\|P_{r_n}A(u(s))-A(u(s))\|^2_Hds\right)^{\frac{1}{2}}\left(\int^t_0\|u^{r_n}(s)\|^2_{H^1}ds\right)^{\frac{1}{2}}\rightarrow 0.
\end{eqnarray*}
On the other hand, the weak convergence of $u^{r_n}$ implies
\begin{eqnarray*}
&&|\int^t_0\langle  A(u(s))(\nabla u^{r_n}(s)-\nabla u(s)), \nabla \phi \rangle ds|\rightarrow
0.
\end{eqnarray*}
Thus, let $n\rightarrow \infty$ in (\ref{equua-1}) to obtain
\begin{eqnarray*}\notag
&&\langle u(t), \phi \rangle-\langle u_0, \phi \rangle-\int^t_0\langle B(u(s)), \nabla \phi \rangle ds\\
\label{equu-5}
&=&-\int^t_0\langle A(u(s))\nabla u(s), \nabla \phi \rangle ds+\int^t_0 \langle\sigma(u(s))h(s), \phi\rangle ds,
\end{eqnarray*}
which means that $u$ is the solution to (\ref{equ-5}). The energy inequality (\ref{equ-6}) is implied by (\ref{equ-09}). To see $u\in C([0,T], H)$, we simply appeal  to Lemma \ref{lem-3}.
\end{flushleft}
$\hfill\blacksquare$

Now, we can define $\mathcal{G}^0: C([0,T],Y)\rightarrow C([0,T],H)\cap L^2([0,T],H^1)$ by
\begin{eqnarray}
\mathcal{G}^0(\check{h}):=\left\{
                   \begin{array}{ll}
                      u_{h}, & {\rm{if}}\ \check{h}= \int^{\cdot}_0 h(s)ds\quad {\rm{for\ some}}\  h\in L^2([0,T],U),\\
                    0, & {\rm{otherwise}}.
                   \end{array}
                  \right.
\end{eqnarray}
\section{Large deviations}
This section is devoted to the proof of the main result.
For any $N>0$ and a given  family $\{h_{\varepsilon}; \varepsilon>0\}\subset \mathcal{A}_N$, let $\bar{u}_{h_{\varepsilon}}$ denote the solution of the following SPDE
\begin{eqnarray}\notag
&&\bar{u}_{h_{\varepsilon}}(t)+\int^t_0div(B(\bar{u}_{h_{\varepsilon}}))ds\\
\label{equ-31}
&&=u_0+\int^t_0div(A(\bar{u}_{h_{\varepsilon}})\nabla\bar{u}_{h_{\varepsilon}})ds
+\sqrt{\varepsilon}\int^t_0\sigma(\bar{u}_{h_{\varepsilon}})dW(s)
 +\int^t_0\sigma(\bar{u}_{h_{\varepsilon}})h_{\varepsilon}(s)ds.
\end{eqnarray}
Then, we have $\bar{u}_{h_{\varepsilon}}(\cdot)=\mathcal{G}^{\varepsilon}(W(\cdot)+\frac{1}{\sqrt{\varepsilon}}\int^{\cdot}_0h_{\varepsilon}(s)ds)
$.

Referring to \cite{DHV}, \cite{H-Z} and using the same method as in Section 4, the following result holds.
\begin{prp}\label{prp-5}
  Assume $u_0\in L^p(\Omega,\mathcal{F}_0, L^p(\mathbb{T}^d))$ for all $p\in [1, \infty)$ and Hypothesis H holds. Then there exists a unique solution $\bar{u}_{h_{\varepsilon}}$ to equation (\ref{equ-31}). Moreover,  we have
\begin{eqnarray*}
\sup_{\varepsilon}\mathbb{E}\Big\{\sup_{0\leq t\leq T}\|\bar{u}_{h_{\varepsilon}}(t)\|^2_H+\int^T_0\|\bar{u}_{h_{\varepsilon}}(t)\|^2_{H^1}dt\Big\}<\infty,\\
\sup_{\varepsilon}\mathbb{E}\sup_{0\leq t\leq T}\|\bar{u}_{h_{\varepsilon}}(t)\|^{2p}_H<\infty, \quad \forall p\geq 1.
\end{eqnarray*}

\end{prp}
The following result concerns the tightness of the above solution.
\begin{prp}\label{prp-3}
The family $\{\bar{u}_{h_{\varepsilon}}, \varepsilon>0\}$ is tight in $L^2([0,T],H).$
\end{prp}
\begin{proof} Set
\[
F(u):=-div(B(u))+div(A(u)\nabla u).
\]
We have
\begin{eqnarray}\notag
\bar{u}_{h_{\varepsilon}}(t)
&=&u_0+\int^t_0F(\bar{u}_{h_{\varepsilon}})ds
+\sqrt{\varepsilon}\int^t_0\sigma(\bar{u}_{h_{\varepsilon}})dW(s) +\int^t_0\sigma(\bar{u}_{h_{\varepsilon}})h_{\varepsilon}(s)ds\\
\label{equ-32}
&:=&J_1+J_2(t)+J_3(t)+J_4(t).
\end{eqnarray}
Clearly, $\mathbb{E}\|J_1\|^2_H\leq C_1$. By (4.2) in \cite{H-Z}, we have
\[
\|F(\bar{u}_{h_{\varepsilon}})\|_{H^{-1}}\leq C(1+\|\bar{u}_{h_{\varepsilon}}\|_{H^1}),
\]
then,
\begin{eqnarray*}
\|J_2(t)-J_2(s)\|_{H^{-1}}
&=&\|\int^t_sF(\bar{u}_{h_{\varepsilon}}(l))dl \|_{H^{-1}}\\
&\leq& \int^t_s\|F(\bar{u}_{h_{\varepsilon}}(l))\|_{H^{-1}}dl\\
&\leq& C\int^t_s(1+\|\bar{u}_{h_{\varepsilon}}\|_{H^1}) dl\\
&\leq& C(t-s)+C(t-s)^{\frac{1}{2}}\Big(\int^t_s\|\bar{u}_{h_{\varepsilon}}\|^2_{H^1}dl\Big)^{\frac{1}{2}}.
\end{eqnarray*}
Thus, combining with Proposition \ref{prp-5}, we have for $\alpha\in (0,\frac{1}{2})$,
\begin{eqnarray*}\notag
\mathbb{E}\|J_2\|_{W^{\alpha,2}([0,T];H^{-1})}
\leq C_2(\alpha).
\end{eqnarray*}
Moreover, by (\ref{equa-3}), we get
\begin{eqnarray*}\notag
&&\mathbb{E}\|J_3(t)-J_3(s)\|^2_H \\ \notag
&=&\varepsilon\mathbb{E} \|\int^t_s\sigma(\bar{u}_{h_{\varepsilon}})dW_r\|^2_H\\ \notag
&\leq&\varepsilon\mathbb{E}\int_{\mathbb{T}^d}\int^t_s\sum_{k\geq 1}|\sigma_k(\bar{u}_{h_{\varepsilon}}(r,z))|^2drdz\\ \notag
&\leq&\varepsilon C\mathbb{E}\int_{\mathbb{T}^d}\int^t_s(1+|\bar{u}_{h_{\varepsilon}}(r,z)|^2)drdz\\
\label{equ-34}
&\leq& C(t-s)\mathbb{E}(1+\sup_{t\in [0,T]}\|\bar{u}_{h_{\varepsilon}}(t)\|^2_H).
\end{eqnarray*}
This implies  that for $\alpha\in (0,\frac{1}{2})$,
\begin{eqnarray*}\notag
\mathbb{E}\|J_3\|_{W^{\alpha,2}([0,T];H)}\leq C_3(\alpha).
\end{eqnarray*}
By H\"{o}lder's inequality and (\ref{equa-3}), we obtain
\begin{eqnarray*}\notag
\|J_4(t)-J_4(s)\|^2_H
&=&\|\int^t_s\sigma(\bar{u}_{h_{\varepsilon}})h_{\varepsilon}(r)dr\|^2_H\\ \notag
&\leq &\int_{\mathbb{T}^d}|\int^t_s\sum_{k\geq 1}|\sigma_k(\bar{u}_{h_{\varepsilon}})||h^k_{\varepsilon}(r)|dr|^2dz\\ \notag
&\leq&\int_{\mathbb{T}^d}|\int^t_s\left(\sum_{k\geq 1}|\sigma_k(\bar{u}_{h_{\varepsilon}})|^2\right)^{\frac{1}{2}}\left(\sum_{k\geq 1}|h^k_{\varepsilon}(r)|^2\right)^{\frac{1}{2}}dr|^2dz\\ \notag
&\leq&\int_{\mathbb{T}^d}\left(\int^t_s\sum_{k\geq 1}|\sigma_k(\bar{u}_{h_{\varepsilon}})|^2dr\right)\left(\int^t_s\sum_{k\geq 1}|h^k_{\varepsilon}(r)|^2dr\right)dz\\ \notag
&\leq & C\int_{\mathbb{T}^d}\int^t_s(1+|\bar{u}_{h_{\varepsilon}}|^2)drdz\int^t_s|h_{\varepsilon}(r)|^2_Udr\\ \label{equ-36}
&\leq &C(t-s)(1+\sup_{t\in [0,T]}\|\bar{u}_{h_{\varepsilon}}(t)\|^2_H)\int^t_s|h_{\varepsilon}(r)|^2_Udr.
\end{eqnarray*}
In view of the boundedness of $h_{\varepsilon}$ in $L^2([0,T], U)$, we deduce that for $\alpha\in (0,\frac{1}{2})$,
\begin{eqnarray*}\label{equ-37}
\mathbb{E}\|J_4\|_{W^{\alpha,2}([0,T];H)}
\leq C_4(\alpha).
\end{eqnarray*}
Collecting the above estimates, we obtain
\begin{eqnarray}\label{equ-38}
\sup_{\varepsilon\in [0,\varepsilon_0]}\mathbb{E}\|\bar{u}_{h_{\varepsilon}}\|_{W^{\alpha,2}([0,T];H^{-1})}<+\infty \quad {\rm{for \ some \ }}  \varepsilon_0>0.
\end{eqnarray}
Set $\Lambda=L^2([0,T],H^1)\cap W^{\alpha,2}([0,T];H^{-1})$. Since the imbedding $H^1(\mathbb{T}^d)\subset H$ is compact, we see from Lemma \ref{lem-1} that $\Lambda$ is compactly  imbedded in $L^2([0,T],H)$.
Set
\[
\|\cdot\|_{\Lambda}:=\|\cdot\|_{L^2([0,T],H^1)}+\|\cdot\|_{W^{\alpha,2}([0,T];H^{-1})}.
\]
Then, for any $L>0$, $K_{L}=\{u\in L^2([0,T],H), \|u\|_{\Lambda}\leq L\}$ is a compact subset of $L^2([0,T],H)$. By Proposition \ref{prp-5} and (\ref{equ-38}), we have
\begin{eqnarray}
\sup_{\varepsilon\in [0,\varepsilon_0]}\mathbb{E}\|\bar{u}_{h_{\varepsilon}}\|_{\Lambda}=C<+\infty.
\end{eqnarray}
Since
\[
\mathbb{P}(\bar{u}_{h_{\varepsilon}}\notin K_L)\leq \mathbb{P}(\|\bar{u}_{h_{\varepsilon}}\|_{\Lambda}\geq L)
\leq \frac{1}{L}\mathbb{E}(\|\bar{u}_{h_{\varepsilon}}\|_{\Lambda})
\leq \frac{C}{L},
\]
 choosing $L$ to be as large as we wish, it follows that $\{\bar{u}_{h_{\varepsilon}}; \varepsilon >0\}$ is tight in $L^2([0,T],H)$.
\end{proof}

Arguing similarly to the above, we have $\{\bar{u}_{h_{\varepsilon}}; \varepsilon >0\}$ are bounded in
 $C([0,T],H)\cap W^{\alpha,p}([0,T];H^{-1})$ for any $p\geq 2$. Choosing $\alpha p>1$, as a consequence, we get
\begin{prp}\label{prp-4}
$\{\bar{u}_{h_{\varepsilon}}; \varepsilon >0\}$ is tight in $C([0,T],H^{-1})$.
\end{prp}
\noindent{\bf Proof of Theorem \ref{thm-7}}.\quad
According to Theorem \ref{thm-6}, the proof of Theorem \ref{thm-7} will be completed if the following Theorem \ref{thm-4} and Theorem \ref{thm-5} are established.

$\hfill\blacksquare$

\begin{thm}\label{thm-4}
Fix $N\in \mathbb{N}$ and let $\{h_\varepsilon\}\subset \mathcal{A}_N$. If $h_\varepsilon$ converge in distribution to $h$ as $\varepsilon\rightarrow 0$, then
\[
\mathcal{G}^{\varepsilon}\left(W(\cdot)+\frac{1}{\sqrt{\varepsilon}}\int^{\cdot}_0h_\varepsilon(s)ds\right)\ converges\ in\ distribution\ to\ \mathcal{G}^{0}\left(\int^{\cdot}_0h(s)ds\right),
\]
in the space $C([0,T],H)$ as $\varepsilon\rightarrow 0$.
\end{thm}
\begin{proof}Recall $\bar{u}_{h_{\varepsilon}}(\cdot)=\mathcal{G}^{\varepsilon}(W(\cdot)+\frac{1}{\sqrt{\varepsilon}}\int^{\cdot}_0h_\varepsilon(s)ds) $. By Proposition \ref{prp-3} and Proposition \ref{prp-4}, we know that $\{\bar{u}_{h_{\varepsilon}}; \varepsilon >0\}$  is tight in the space
$L^2([0,T],H)\cap C([0,T],H^{-1})$. Let $(u,h, W)$ be any limit point of the tight family $\{(\bar{u}_{h_{\varepsilon}},h_{\varepsilon},W), \varepsilon\in(0,\varepsilon_0)\}$, we have to show that $u$ has the same law as $u_h=\mathcal{G}^{0}(\int^{\cdot}_0h(s)ds)$, and $\bar{u}_{h_{\varepsilon}}\rightarrow u$ in distribution in the space $C([0,T],H)$.

Set
\[
\Pi=\Big(L^2([0,T],H)\cap C([0,T],H^{-1}), S_N, C([0,T],Y)\Big).
\]
By the Skorokhod representation theorem, there exists a stochastic basis $(\Omega^1, \mathcal{F}^1, \{\mathcal{F}^1_t\}_{t\in [0,T]}, \mathbb{P}^1)$ and $\Pi-$valued random variables $(\tilde{X}_{\varepsilon}, \tilde{h}_{\varepsilon},\tilde{W}_{\varepsilon}),(\tilde{X}, \tilde{h},\tilde{W})$ on this basis, such that $(\tilde{X}_{\varepsilon}, \tilde{h}_{\varepsilon},\tilde{W}_{\varepsilon})$ (resp. $(\tilde{X}, \tilde{h},\tilde{W})$) has the same law as $(\bar{u}_{h_{\varepsilon}},h_{\varepsilon},W)$ (resp. $(u,h,W)$), and $(\tilde{X}_{\varepsilon}, \tilde{h}_{\varepsilon},\tilde{W}_{\varepsilon})\rightarrow (\tilde{X}, \tilde{h},\tilde{W})$, $\mathbb{P}^1-$a.s. in $\Pi$. Hence,
\begin{eqnarray}\label{e-71}
\int^T_0 \|\tilde{X}_{\varepsilon}-\tilde{X}\|^2_Hdt\rightarrow 0\quad  \mathbb{P}^1-a.s..
\end{eqnarray}
Denote by $\mathbb{E}^1$ the expectation associated with the stochastic basis $(\Omega^1, \mathcal{F}^1, \{\mathcal{F}^1_t\}_{t\in [0,T]}, \mathbb{P}^1)$.  In the following, to simplify the notations, we denote
 \[
X=\tilde{X},\ h= \tilde{h},\ W=\tilde{W},\ {X}_{\varepsilon}=\tilde{X}_{\varepsilon}, \ h_\varepsilon=\tilde{h}_{\varepsilon},\ {W}_{\varepsilon}=\tilde{W}_{\varepsilon}.
\]
From the equation satisfied by $(\bar{u}_{h_{\varepsilon}},h_{\varepsilon},W)$, it follows that $({X}_{\varepsilon},{h}_{\varepsilon}, {W}_{\varepsilon})$ satisfies the following equation
\begin{eqnarray}\notag
{X}_{\varepsilon}(t)
&=&u_0-\int^t_0div(B({X}_{\varepsilon}))ds+\int^t_0div(A({X}_{\varepsilon})\nabla{X}_{\varepsilon})ds\\
\label{equ-39}
&&+\sqrt{\varepsilon}\int^t_0\sigma({X}_{\varepsilon})dW_{\varepsilon}(s) +\int^t_0\sigma({X}_{\varepsilon})h_{\varepsilon}(s)ds.
\end{eqnarray}
Similar to the proof of Theorem \ref{thm-2}, using (\ref{e-71}), let $\varepsilon \rightarrow 0$ in (\ref{equ-39}) to see that $X$ is the unique solution of the following equation
\begin{eqnarray}\label{equ-40}
{X}(t)=u_0-\int^t_0div(B({X}))ds+\int^t_0div(A({X})\nabla{X})ds
+\int^t_0\sigma({X})h(s)ds.
\end{eqnarray}
By the uniqueness of the skeleton equation, we have $X=u_h$, which implies that $u$ has the same law as $\mathcal{G}^{0}(\int^{\cdot}_0h(s)ds)$. Moreover, we have $X\in C([0,T],H)\cap L^2([0,T],H^1)$.

Finally, to complete the proof, it suffices to show that
\[
\lim_{\varepsilon\rightarrow 0}\sup_{t\in [0,T]}\|X_{\varepsilon}(t)-X(t)\|_H=0, \quad {\rm{in\ probability}}.
\]
Let $\omega_{\varepsilon}(t)=X_{\varepsilon}(t)-X(t)$. We have
\begin{eqnarray}\notag
&&d\omega_{\varepsilon}(t)
+(div(B({X}_{\varepsilon}))-div(B({X})))dt\\
\label{equ-41}
&=&(div(A({X}_{\varepsilon})\nabla{X}_{\varepsilon})-div(A({X})\nabla{X}))dt
+\sqrt{\varepsilon}\sigma({X}_{\varepsilon})dW_{\varepsilon}(t) +(\sigma({X}_{\varepsilon})h_{\varepsilon}-\sigma(X)h)dt.
\end{eqnarray}
Applying It\^{o} formula, we obtain
\begin{eqnarray}\notag
\|\omega_{\varepsilon}(t)\|^2_H&=&2\int^t_0\langle\omega_{\varepsilon},-div(B({X}_{\varepsilon})-B({X}))\rangle ds\\ \notag
&&+2\int^t_0\langle\omega_{\varepsilon},div(A({X}_{\varepsilon})\nabla{X}_{\varepsilon}-A({X})\nabla{X})\rangle ds\\ \notag
&&+2\int^t_0\langle\omega_{\varepsilon}, \sigma({X}_{\varepsilon})h_{\varepsilon}-\sigma(X)h\rangle ds+2\int^t_0\langle\omega_{\varepsilon},\sqrt{\varepsilon}\sigma({X}_{\varepsilon})dW_{\varepsilon}(s)\rangle \\ \notag
&&+\varepsilon \int^t_0 \int_{\mathbb{T}^d}\sum^{\infty}_{k=1}|\sigma_k({X}_{\varepsilon}(s,z))|^2dzds,\\ \notag
&=&2\int^t_0\langle \nabla \omega_{\varepsilon},B({X}_{\varepsilon})-B({X})\rangle ds-2\int^t_0\langle \nabla \omega_{\varepsilon},A({X}_{\varepsilon})\nabla{X}_{\varepsilon}-A({X})\nabla{X}\rangle ds\\ \notag
&&+2\int^t_0\langle\omega_{\varepsilon}, \sigma({X}_{\varepsilon})h_{\varepsilon}-\sigma(X)h\rangle ds
+2\int^t_0\langle\omega_{\varepsilon},\sqrt{\varepsilon}\sigma({X}_{\varepsilon})dW_{\varepsilon}(s)\rangle \\ \notag
&&+\varepsilon \int^t_0 \int_{\mathbb{T}^d}\sum^{\infty}_{k=1}|\sigma_k({X}_{\varepsilon}(s,z))|^2dzds, \\ \label{equ-42}
&:=& K_1(t)+K_2(t)+K_3(t)+K_4(t)+K_5(t).
\end{eqnarray}
By Young's inequality, for $\delta>0$, we have
\begin{eqnarray}\notag
K_1(t) &\leq& C\int^t_0\|\omega_{\varepsilon}\|_H\|\omega_{\varepsilon}\|_{H^1}ds\\
\label{equ-43}
&\leq& \frac{\delta}{4}\int^t_0\|\omega_{\varepsilon}\|^2_{H^1}ds+C\int^t_0\|\omega_{\varepsilon}\|^2_Hds.
\end{eqnarray}
 $K_2(t)$ can be written as
\begin{eqnarray*}\notag
K_2(t)
&=&-2\int^t_0\langle \nabla \omega_{\varepsilon}, A({X}_{\varepsilon})\nabla{X}_{\varepsilon}-A({X}_{\varepsilon})\nabla X\rangle ds\\ \notag
&&-2\int^t_0\langle \nabla \omega_{\varepsilon}, A({X}_{\varepsilon})\nabla X-A({X})\nabla X\rangle ds\\ \label{equ-44}
&:=& L_1(t)+L_2(t).
\end{eqnarray*}
In view of  Hypothesis H, we have
\begin{eqnarray}\label{equ-48}
L_1(t)\leq -2\delta \int^t_0\|\omega_{\varepsilon}\|^2_{H^1}ds.
\end{eqnarray}
By H\"{o}lder's inequality and Young's inequality, we have
\begin{eqnarray}\notag
L_2(t)&\leq& \frac{\delta}{4}\int^t_0\|\omega_{\varepsilon}\|^2_{H^1}ds+C\int^t_0\int_{\mathbb{T}^d}|A({X}_{\varepsilon})-A({X})|^2|\nabla X|^2dzds\\
\label{equ-47}
&:=& \frac{\delta}{4}\int^t_0\|\omega_{\varepsilon}\|^2_{H^1}ds+N^{\varepsilon}(t).
\end{eqnarray}
 For any constant $M>0$, we have
 \begin{eqnarray}\notag
N^{\varepsilon}(t)&=& C\int^t_0\int_{\mathbb{T}^d}|A({X}_{\varepsilon})-A({X})|^2|\nabla X|^2I_{\{|\nabla X(z)|\leq M\}}dzds\\ \notag
&&\ +C\int^t_0\int_{\mathbb{T}^d}|A({X}_{\varepsilon})-A({X})|^2|\nabla X|^2I_{\{|\nabla X(z)|>M\}}dzds\\ \notag
&\leq & CM^2\int^t_0\int_{\mathbb{T}^d}|A({X}_{\varepsilon})-A({X})|^2dzds+C\int^t_0\int_{\mathbb{T}^d}|\nabla X(z)|^2I_{\{|\nabla X(z)|>M\}}dzds\\
\label{equat-1}
&\leq& CM^2\int^t_0\|\omega_{\varepsilon}\|^2_Hds+C\int^t_0\int_{\mathbb{T}^d}|\nabla X(z)|^2I_{\{|\nabla X(z)|>M\}}dzds.
\end{eqnarray}
Since $X\in L^2([0,T],H^1)$, it yields
\[
\lim_{M\rightarrow \infty}\int^t_0\int_{\mathbb{T}^d}|\nabla X(z)|^2I_{\{|\nabla X(z)|>M\}}dzds=0.
\]
In light of (\ref{e-71}), it follows that
\begin{equation}\label{2}
\lim_{\varepsilon\rightarrow 0}N^{\varepsilon}(T)=0, \quad \mathbb{P}-a.s..
\end{equation}
Combining (\ref{equ-48}) and (\ref{2}), we have
\begin{eqnarray}\label{equu-6}
K_2(t)\leq-\frac{7\delta}{4} \int^t_0\|\omega_{\varepsilon}\|^2_{H^1}ds+N^{\varepsilon}(T),\quad {\rm{and}} \quad N^{\varepsilon}(T) \rightarrow 0, \quad \mathbb{P}-a.s..
\end{eqnarray}
 By H\"{o}lder's inequality and (\ref{equa-3}), we bound $K_3(t)$ in (\ref{equ-42}) as follows
\begin{eqnarray*}
K_3(t)&=& 2\int^t_0\langle\omega_{\varepsilon}, \sigma({X})(h_{\varepsilon}-h)\rangle ds+2\int^t_0\langle\omega_{\varepsilon},(\sigma({X}_{\varepsilon})-\sigma(X))h_{\varepsilon}\rangle ds \\ \notag
&\leq&2\int^t_0|h_{\varepsilon}-h|_U\|\omega_{\varepsilon}\|_H\left(\int_{\mathbb{T}^d}\sum_{k\geq 1}|\sigma_k({X})|^2dz\right)^{\frac{1}{2}}ds\\
&&\ +2\int^t_0|h_{\varepsilon}|_U\|\omega_{\varepsilon}\|_H\left(\int_{\mathbb{T}^d}\sum_{k\geq 1}|\sigma_k({X}_{\varepsilon})-\sigma_k({X})|^2dz\right)^{\frac{1}{2}}ds\\ \notag
&\leq& C\int^t_0|h_{\varepsilon}-h|_U\|\omega_{\varepsilon}\|_H(1+\|{X}(s)\|_H)ds+C\int^t_0|h_{\varepsilon}|_U\|\omega_{\varepsilon}\|^2_Hds\\ \notag
&\leq&  C\sup_{t\in [0,T]}(1+\|{X}\|_H)\int^t_0\|\omega_{\varepsilon}\|_H|h_{\varepsilon}-h|_Uds+C\int^t_0\|\omega_{\varepsilon}\|^2_H |h_{\varepsilon}|_Uds\\ \notag
&\leq&  C\sup_{t\in [0,T]}(1+\|{X}\|_H)\left(\int^t_0|h_{\varepsilon}-h|^2_Uds\right)^{\frac{1}{2}}\left(\int^t_0\|\omega_{\varepsilon}\|^2_Hds\right)^{\frac{1}{2}}+2C\int^t_0\|\omega_{\varepsilon}\|^2_H |h_{\varepsilon}|_Uds\\
&\leq&  C N^{\frac{1}{2}}\left(\int^t_0\|\omega_{\varepsilon}\|^2_Hds\right)^{\frac{1}{2}}+2C\int^t_0\|\omega_{\varepsilon}\|^2_H |h_{\varepsilon}|_Uds.
\end{eqnarray*}
By the Burkholder-Davis-Gundy inequality and (\ref{equa-3}), we have
\begin{eqnarray*}
\mathbb{E}\sup_{t\in [0,T]}|K_4(t)|
&\leq& 2\sqrt{\varepsilon}C\mathbb{E}\left[\sum_{k\geq 1}\int^T_0\left(\int_{\mathbb{T}^d}\omega_{\varepsilon}\sigma_k({X}_{\varepsilon}(s,z))dz\right)^2ds\right]^{\frac{1}{2}}\\
&\leq&2\sqrt{\varepsilon}C\mathbb{E}\left[\sum_{k\geq 1}\int^T_0\left(\int_{\mathbb{T}^d}|\omega_{\varepsilon}|^2dz\right)\left(\int_{\mathbb{T}^d}|\sigma_k({X}_{\varepsilon}(s,z))|^2dz\right)ds\right]^{\frac{1}{2}}\\
&\leq&2\sqrt{\varepsilon}C\mathbb{E}\left[\int^T_0\|\omega_{\varepsilon}\|^2_H(1+\|{X}_{\varepsilon}(s)\|^2_H)ds\right]^{\frac{1}{2}}\\
&\leq&2\sqrt{\varepsilon}C\mathbb{E}\sup_{t\in[0,T]}(1+\|{X}_{\varepsilon}(t)\|_H) \left(\int^T_0\|\omega_{\varepsilon}(s)\|^2_Hds\right)^{\frac{1}{2}}\\
&\leq&\sqrt{\varepsilon}C\Big(\mathbb{E}\sup_{t\in[0,T]}(1+\|{X}_{\varepsilon}(t)\|^2_H)\Big)^{\frac{1}{2}} \left(\mathbb{E}\int^T_0(\|X_{\varepsilon}(s)\|^2_H+\|X(s)\|^2_H)ds\right)^{\frac{1}{2}}\\
&\leq&\sqrt{\varepsilon}C,
\end{eqnarray*}
which implies
\begin{eqnarray}\label{eq-6}
\sup_{t\in [0,T]}|K_4(t)|\rightarrow 0,
\end{eqnarray}
in probability as $\varepsilon\rightarrow 0$.

Utilizing the linear growth of $\sigma$, we have
\begin{eqnarray*}
K_5(t)&=&\varepsilon \int^t_0 \int_{\mathbb{T}^d}\sum^{\infty}_{k=1}|\sigma_k({X}_{\varepsilon}(s,z))|^2dzds\\
&\leq& \varepsilon C\int^t_0 \int_{\mathbb{T}^d}(1+|{X}_{\varepsilon}(s,z)|^2) dzds\\ \notag
&\leq& \varepsilon C\int^t_0 (1+\|{X}_{\varepsilon}(s)\|^2_H)ds,
\end{eqnarray*}
then, it yields
\begin{eqnarray}\label{eq-8}
K_5(T)\rightarrow 0,
\end{eqnarray}
in probability as $\varepsilon\rightarrow 0$.

Collecting all the above estimates, we obtain
\begin{eqnarray}\notag
\|\omega_{\varepsilon}(t)\|^2_H
&\leq& C\int^t_0\|\omega_{\varepsilon}\|^2_Hds+N^{\varepsilon}(T)+C N^{\frac{1}{2}}\left(\int^T_0\|\omega_{\varepsilon}\|^2_Hds\right)^{\frac{1}{2}}\\  \notag
&&+2C\int^t_0\|\omega_{\varepsilon}\|^2_H |h_{\varepsilon}|_Uds
+\sup_{t\in [0,T]}|K_4(t)|+K_5(T) \\
\label{equ-45}
&:=& C\int^t_0(1+|h_{\varepsilon}|_U)\|\omega_{\varepsilon}\|^2_Hds+\Theta(\varepsilon, T),
\end{eqnarray}
where
\begin{eqnarray*}
\Theta(\varepsilon, T)= N^{\varepsilon}(T)+C N^{\frac{1}{2}}\left(\int^T_0\|\omega_{\varepsilon}\|^2_Hds\right)^{\frac{1}{2}} +\sup_{t\in [0,T]}|K_4(t)|+K_5(T).
\end{eqnarray*}
Applying Gronwall's inequality, it follows from (\ref{equ-45}) that
\begin{eqnarray*}
\sup_{t\in [0,T]}\|\omega_{\varepsilon}(t)\|^2_H
\leq \Theta(\varepsilon, T)\exp\Big\{C\int^T_0(1+|h_{\varepsilon}(s)|_U)ds\Big\}.
\end{eqnarray*}
By (\ref{e-71}), (\ref{2}), (\ref{eq-6}) and (\ref{eq-8}), we have
\[
\Theta(\varepsilon, T)\rightarrow0,\ \ {\rm{in \ probability}} \quad {\rm{as}}\ \ \varepsilon\rightarrow 0.
\]
Since $\int^T_0 |h_{\varepsilon}(s)|^2_Uds\leq N, \mathbb{P}-a.s.$, then
\[
\exp\Big\{\int^T_0C(1+|h_{\varepsilon}(s)|_U)ds\Big\}<\infty,\quad  \mathbb{P}-a.s.
\]
we conclude that
\[
\sup_{t\in [0,T]}\|\omega_{\varepsilon}(t)\|^2_H\rightarrow 0,
\]
in probability as $\varepsilon\rightarrow 0$. We complete the proof.
\end{proof}

Replacing $\sqrt{\varepsilon}\int^t_0\sigma({X}_{\varepsilon})dW_{\varepsilon}$ by 0 in the proof of Proposition \ref{prp-3}, Proposition \ref{prp-4} and Theorem \ref{thm-4}, we can prove the following result.
\begin{thm}\label{thm-5}
$\mathcal{G}^0(\int^{\cdot}_0h(s)ds)$ is a continuous mapping from $h\in S_N$ into $C([0,T],H)$, in particular, $\{\mathcal{G}^0(\int^{\cdot}_0h(s)ds); h\in S_N\}$ is a compact subset of $C([0,T],H)$.
\end{thm}


\

\noindent{\bf  Acknowledgements}\  The authors are grateful to the anonymous referees for comments and suggestions. This work is partly supported by National Natural Science Foundation of China (No.11371041, 11671372, 11431014, 11401557, 11801032). Key Laboratory of Random Complex Structures and Data Science, Academy of Mathematics and Systems Science, Chinese Academy of Sciences (No. 2008DP173182). China Postdoctoral Science Foundation funded project (No. 2018M641204).

\def\refname{ References}

\end{document}